\documentclass{gtart_h}

\def\ifplaintex{\expandafter\ifx\csname documentclass\endcsname\relax}

\def\gtp{{\mathsurround=0pt\it $\cal G\mskip-2mu$eometry \&\ 
$\cal T\!\!$opology $\cal P\!$ublications}}  

\def\recd{{\small Received:\qua\receiveddate\ifx\reviseddate\relax
\else\qquad Revised:\qua\reviseddate\fi\par}} 

\def\lognumber#1{\def\thelognumber{#1}}
\def\volumenumber#1{\def\thevolumenumber{#1}}
\def\volumeyear#1{\def\thevolumeyear{#1}}
\def\papernumber#1{\def\thepapernumber{#1}}
\def\pagenumbers#1#2{\def\startpage{#1}\def\finishpage{#2}}
\def\published#1{\def\publishdate{#1}}

\def\received#1{\def\receiveddate{#1}}
\def\revised#1{\def\reviseddate{#1}}
\def\accepted#1{\def\accepteddate{#1}}

\let\\\par\let\thelognumber\relax\let\thevolumenumber\relax
\let\thepapernumber\relax\let\thevolumeyear\relax\let\startpage\relax
\let\finishpage\relax\let\publishdate\relax\let\receiveddate\relax
\let\reviseddate\relax\let\accepteddate\relax\let\theasciititle\relax
\let\theasciiauthors\relax
\let\theasciiabstract\relax

\let\theasciiemail\relax

\ifplaintex
\font\logobig=cmssbx10 scaled 3836
\font\logomed=cmssbx10 scaled 2557
\else
\font\logobig=cmssbx10 scaled 4200
\font\logomed=cmssbx10 scaled 2800
\fi

\long\def\makeagttitle{   
\count0=\startpage
\agt\hfill      
\hbox to 45truept{\vbox to 0pt{\vglue -13truept{\logomed A\kern -.37em{\logobig 
T}\kern -.38em G}\vss}\hss}
\break
{\small Volume \thevolumenumber\ (\thevolumeyear)
\startpage--\finishpage\nl
Published: \publishdate}

\vglue .25truein

{\parskip=0pt\leftskip 0pt plus
1fil\def\\{\par\smallskip}{\Large\bf\thetitle}\par\medskip} \vglue
0.05truein

{\parskip=0pt\leftskip 0pt plus 1fil\def\\{\par}{\sc\theauthors}
\par\medskip}%
 
\vglue 0.03truein 

{\small\leftskip 25truept\rightskip 25truept{\bf Abstract}\stdspace\theabstract

{\bf AMS Classification}\stdspace\theprimaryclass
\ifx\thesecondaryclass\relax\else; \thesecondaryclass\fi\par
{\bf Keywords}\stdspace \thekeywords\par}\vglue 7truept

}   

\ifplaintex
\font\phead=cmsl9 scaled 950
\font\pnum=cmbx10 scaled 913
\font\pfoot=cmsl9 scaled 950
\headline{\vbox to 0pt{\vskip -4.5mm\line{\small\phead\ifnum
\count0=\startpage ISSN 1472-2739 (on-line) 1472-2747 (printed)
\hfill {\pnum\folio}\else\ifodd\count0\def\\{ }%
\ifx\theshorttitle\relax\thetitle\else\theshorttitle\fi\hfill{\pnum\folio}
\else\def\\{ and }{\pnum\folio}\hfill\ifx\theshortauthors\relax\theauthors
\else\theshortauthors\fi\fi\fi}\vss}}
\footline{\vbox to 0pt{\vglue 0mm\line{\small\pfoot\ifnum\count0=\startpage
\copyright\ \gtp\hfill\else
\agt, Volume \thevolumenumber\ (\thevolumeyear)\hfill\fi}\vss}}
\else
\font=cmsl9 scaled 1050
\font\lnum=cmbx10 
\font=cmsl9 scaled 1050
\makeatletter
\def\@oddhead{{\small\lhead\ifnum\count0=\startpage ISSN 1472-2739 
(on-line) 1472-2747 (printed)\hfill {\lnum\number\count0}\else\ifodd\count0
\def\\{ }\ifx\theshorttitle\relax \thetitle \else\theshorttitle\fi\hfill
{\lnum\number\count0}\else\def\\{ and }{\lnum\number\count0}
\hfill\ifx\theshortauthors\relax 
\theauthors\else\theshortauthors\fi\fi\fi}}\def\@evenhead{\@oddhead}
\def\@oddfoot{\small\lfoot\ifnum\count0=\startpage\copyright\ \gtp\hfill\else
\agt, Volume \thevolumenumber\ (\thevolumeyear)\hfill\fi}
\def\@evenfoot{\@oddfoot}
\makeatother
\fi
\let\maketitlepage\makeagttitle
\let\makeshorttitle\maketitlepage
\let\maketitle\maketitlepage

\newwrite\gtoutfile
\long\gdef\makeheadfile{  
{\def\\{, }\def\s{ }
\immediate\openout\gtoutfile head.xxx
\immediate\write\gtoutfile{Proxy-for: \ifx\theasciiauthors\relax
\theauthors\else\theasciiauthors\fi\s<\ifx\theasciiemail\relax\theemail\else\theasciiemail\fi>}
\immediate\write\gtoutfile{\noexpand\\}
\immediate\write\gtoutfile{Authors: \ifx\theasciiauthors\relax
\theauthors\else\theasciiauthors\fi}
{\def\\{ }\immediate\write\gtoutfile{Title: \ifx\theasciititle\relax
\thetitle\else\theasciititle\fi}}
\immediate\write\gtoutfile{Subj-class: GT or SG, GR etc}
\immediate\write\gtoutfile{MSC-class: \theprimaryclass\ifx\thesecondaryclass\relax\else, \thesecondaryclass\fi}
\immediate\write\gtoutfile{Journal-ref: Algebr. Geom. Topol. \thevolumenumber\s
(\thevolumeyear) \startpage-\finishpage}
\immediate\write\gtoutfile{Comments: Published by Algebraic and
Geometric Topology at}
\immediate\write\gtoutfile{\s\s\s  http://www.maths.warwick.ac.uk/agt/AGTVol\thevolumenumber/agt-\thevolumenumber-\thepapernumber.abs.html}
\immediate\write\gtoutfile{\noexpand\\}
\immediate\write\gtoutfile{}
\ifx\theasciiabstract\relax
\immediate\write\gtoutfile{\theabstract}\else
\immediate\write\gtoutfile{\theasciiabstract}\fi
\immediate\write\gtoutfile{}
\immediate\write\gtoutfile{\noexpand\\}
\immediate\write\gtoutfile{}
\immediate\closeout\gtoutfile}}  

\def\maketitlepage{\makeagttitle\makeheadfile}
\let\makeshorttitle\maketitlepage
\let\maketitle\maketitlepage

\lognumber{63}
\volumenumber{5}
\volumeyear{2005}
\papernumber{63}
\pagenumbers{1573}{1584}
\received{27 January 2005} 
\revised{9 June 2005}
\accepted{8 November 2005}
\published{24 November 2005}

\usepackage{amsmath, amssymb}
\usepackage{graphicx, psfrag}

\def\psfraga <#1,#2> #3#4{%
\psfrag {#3}{\smash{\rlap{\kern #1 \raise #2\hbox{#4}}}}}

\def\figref#1{\hyperlink{#1anchor}{Figure~\ref*{#1}}}
\def\anchor#1{\noindent\hypertarget{#1anchor}{\smash{$\phantom{99}$}}}

\newtheorem{Thm}{Theorem}[section]
\newtheorem*{MThm}{Main Theorem}
\newtheorem{Def}[Thm]{Definition}
\newtheorem{Coro}[Thm]{Corollary}
\newtheorem{Lem}[Thm]{Lemma}
\newtheorem{Prop}[Thm]{Proposition}

\begin{document}

\title{Locally unknotted spines of Heegaard splittings}
\author{Jesse Johnson}
\address{Mathematics Department, University of 
California\\Davis, CA 95616, USA}
\email{jjohnson@math.ucdavis.edu}

\begin{abstract}
We show that under reasonable conditions, the spines of the handlebodies
of a strongly irreducible Heegaard splitting will intersect a closed ball
in a graph which is isotopic into the boundary of the ball.  This is in some 
sense a generalization of the
results by Scharlemann on how a strongly irreducible Heegaard splitting
surface can intersect a ball.
\end{abstract}


\maketitle

\section{Introduction}
\label{defsect}
The goal of this paper is to understand how a strongly irreducible Heegaard
splitting of a manifold $M$ intersects an arbitrary closed ball in
$B \subset M$. Intuitively, a strongly irreducible splitting is as
efficient as possible so the intersection should be very simple.

This problem was first studied by Scharlemann~\cite{schar:local},
who showed that under certain conditions, the splitting surface will intersect
$B$ as the boundary of a neighborhood of a graph on $n$ points which is
isotopic into $\partial B$.

We will look at how the spine of a Heegaard splitting intersects $B$.
In particular, we will show that under reasonable conditions,
the pieces of the spine inside the ball will be unknotted, i.e.\ simultaneously
isotopic into the boundary of the ball.

We will begin by generalizing the notion of a tangle.  For a graph $T$,
let $\partial T$ be the set of valence-one vertices.

\begin{Def}\textup{
Let $B$ be a closed ball.  A \textit{graph tangle} is a graph $T$ (not 
necessarily connected) properly embedded in $B$, i.e.\ such that any 
valence-one vertices are in the boundary of $B$.
}\end{Def}

\begin{Def}\textup{
A graph tangle $T$ is \textit{unknotted} if there is an embedding of
$T \times [0,1]$ such that $T \times \{0\}$ maps to $T$ and
$(T \times \{1\}) \cup (\partial T \times [0,1])$ is embedded in $\partial B$.
In other words, an unknotted graph tangle can be isotoped into the boundary of
$B$.  
}\end{Def}

\begin{Def}\textup{
Let $G$ be a graph and let $P \subset G$ be a finite collection of points,
disjoint from the vertices of $G$, which separate $G$.  A component of
$G \setminus P$ will be called a \textit{cut graph} of $G$.
}\end{Def}

Let $M$ be a compact 3-manifold. Let $B$ be a ball embedded in $M$
and let $K$ be a graph embedded in $M$ such that $B \cap K$ is
non-empty.

\begin{Def}\textup{
A \textit{trivial cone} of $K$ with respect to $B$ is a cut graph $C$ of $K$ 
(usually a cone on
two or more points) and an embedding of $C \times I$
such that $C \times \{0\} = C$,
$(C \times \{1\}) \cup (\partial C \times I) \subset \partial B$,
the interior of $C \times I$ is disjoint from $K$
and near $(C \times \{1\}) \cup (\partial C \times I)$, 
$C \times I$ lies outside of $B$.  (Here, $I = [0,1]$.)
}\end{Def}

\begin{Def}\textup{
A graph, $K$, (not necessarily connected) in a manifold, $M$, is
\textit{locally unknotted} if for any closed ball $B$ such that there are
no trivial cones with respect to $B$, $K \cap B$ is unknotted in $B$.  
Otherwise, $K$ is \textit{locally knotted}.
}\end{Def}

The requirement that there be no trivial cones is necessary because of
the following example:
\begin{figure}[ht!]\anchor{locknot}\small
  \begin{center}
\psfrag {trivial cone}{trivial cone}  
\psfrag {B}{$B$}  
\psfrag {B'}{$B'$}  
\includegraphics[width=2in]{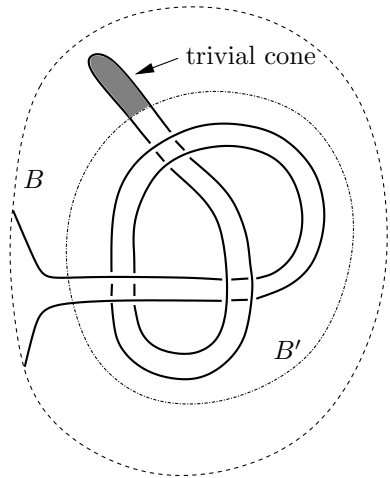}
  \caption{A graph which is knotted in a ball because of a trivial cone}
  \label{locknot}
  \end{center}
\end{figure}
Let $K$ be a graph embedded in $M$.
Let $B$ be a small ball which intersects an edge of $K$ in an unknotted
arc $\alpha$.  Isotope $\alpha$ in $B$ to form the shape shown in
\figref{locknot}.  Let $B'$ be the smaller ball, indicated in the
Figure.  There are two parallel sub-arcs of $\alpha$ in $B'$, each of which is
knotted.  The graph $K$ may still be locally unknotted because
for this choice of $B'$, there is a trivial cone.

\begin{Def}\textup{
Let $M$ be a compact, connected, closed, orientable 3-manifold.
Let $\Sigma$ be a compact, connected, closed, orientable,
two-sided surface in $M$ such that $M \setminus \Sigma$ consists
of two components, $X_1$ and $X_2$.  Assume $H_1 = X_1 \cup \Sigma$ and
$H_2 = X_2 \cup \Sigma$ are handlebodies. Then the ordered triple
$(\Sigma, H_1, H_2)$ is a \textit{Heegaard splitting} for $M$.
}\end{Def}

The splitting surface $\Sigma$ defines the handlebodies $H_1$ and $H_2$ up to labeling so we will often refer to the Heegaard surface $\Sigma$ to indicate the Heegaard splitting $(\Sigma, H_1, H_2)$.

A Heegaard splitting is \textit{stabilized} if there are properly embedded,
essential disks $D_1$, $D_2$ in $H_1$, $H_2$ respectively such that
$\partial D_1$ and $\partial D_2$ intersect transversely in a single point in 
$\Sigma$.  
A splitting is \textit{reducible} (\textit{weakly reducible}) if there are 
properly embedded, essential disks in $H_1$, $H_2$, respectively, whose 
boundaries coincide (are disjoint.)

A \textit{spine} for a handlebody $H_1$ is a graph $K_1$ embedded in $H_1$
such that $H_1 \setminus K_1 \cong \partial H_1 \times (0,1]$.
For a Heegaard splitting $(\Sigma, H_1, H_2)$, a spine for $\Sigma$
is the disjoint union $K = K_1  \cup K_2$ where $K_1$ is a spine for
$H_1$ and $K_2$ is a spine for $H_2$.

We will prove the following:

\begin{MThm}
Let $(\Sigma, H_1, H_2)$ be a strongly irreducible Heegaard splitting for a
manifold $M$.  Then every spine of $\Sigma$ is locally unknotted.
\end{MThm}

In Section~\ref{sweepsect}, we define sweep-outs and introduce the
Rubinstein-Scharlemann graphic.  This is the main tool used in
Section~\ref{lemsect} to prove Proposition~\ref{mainlem}.
The Proposition states that if a Heegaard splitting is strongly irreducible then some level surface of any sweep-out intersects a given ball in trivial pieces.
This is essentially a modification of the result in \cite{rub:compar} 
that given any two sweep-outs, there is a pair of level surfaces such that 
each loop of intersection is either trivial in both surfaces or non-trivial 
in both surfaces.  In our case, $\partial B$ contains only trivial loops so
the loops must all be trivial in the level surface of $f$ as well.
The proof of the Main Theorem is given in the final section.

To indicate the significance of this theorem, we note that if $\Sigma$ is
stabilized then there is a spine of $\Sigma$ which is locally knotted.
\begin{figure}[ht!]\anchor{stabknot}\small
  \begin{center}
\psfraga <-3pt,0pt> {B}{$B$}  
\psfrag {e}{$e$}  
\psfrag {e'}{$e'$}  
\includegraphics[width=2in]{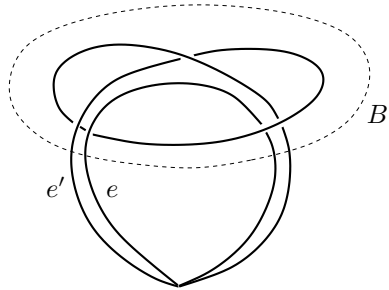}
  \caption{A locally knotted spine of a stabilized Heegaard splitting}
  \label{stabknot}
  \end{center}
\end{figure}
To see this, let $(\Sigma, H_1, H_2)$ be a stabilization of a strongly irreducible
Heegaard splitting $(\Sigma', H_1', H_2')$.
Let $K$ be a spine of $H_1'$ with a single vertex. Because $\Sigma'$ is 
strongly irreducible, any edge $e$ of $K$ is not isotopic into a ball.  
Form a spine for $H_1$ by adding a trivial loop $e'$
to $K$.  Slide $e'$ along $e$ to create a formation as in \figref{stabknot}
and let $B$ be the ball indicated by the dotted line.  The two arcs are
knotted in $B$ and there are no trivial cones (because $e$ is not isotopic
into a ball) so $H_1$ is locally knotted.

This research was supported by NSF {\small VIGRE} grant 0135345.

\section{The Rubinstein-Scharlemann graphic}
\label{sweepsect}

Let $M$ be a compact, connected, closed, orientable 3-manifold and $(\Sigma, H_1, H_2)$ a
Heegaard splitting for $M$.  Let $K_i$ be a spine for $H_i$.  Then
$H_i \setminus K_i$ is homeomorphic to $\Sigma \times (0,1]$.
Let $h : (H_i \setminus K_i) \rightarrow \Sigma \times (0,1] :
         x \mapsto (p_x, l_x)$ be a homeomorphism and let
$f : H_i \rightarrow I$ such that $f(x) = l_x$ when $x \notin K_i$ and
$f(x) = 0$ when $x \in K_i$.  If we combine the two functions (one for each
handlebody,) we get a function on all of $M$.

\begin{Def}\textup{
A \textit{sweep-out} is a smooth function $f : M \rightarrow [-1,1]$ such that
$f^{-1}(-1)$ is a spine of $H_1$, $f^{-1}(1)$ is a spine of $H_2$ and
$f^{-1}(t)$ is a surface isotopic to $\Sigma$ for each $t \in (-1,1)$.
}\end{Def}

By choosing a smooth embedding for each $K_i$ and a smooth map $h$, the above construction produces a sweep-out for any Heegaard splitting of $M$.

\begin{Def}\textup{
Let $f$ and $g$ be sweep-outs and define $f_t = f^{-1}(t)$ and 
$g_t = g|_{f_t}$.  We will say that $f$ and $g$ are in
\textit{general position} if $g_t$ is a Morse function on $f_t$ for
all but finitely many $t$, and $g_t$ is a near-Morse function for the remaining
values.  (A function is \textit{near-Morse} if there is a
single degenerate critical point or there are two non-degenerate critical
points at the same level.)
}\end{Def}

Note that for $t \in (-1,1)$, $f_t$ is a surface embedded in $M$, while $g_t$ is a function on $f_t$.

Let $f$ and $g$ be sweep-outs in general position.  For each $t \in [-1,1]$,
let $s_t$ be a critical value of the function $g_t$.  Rubinstein and
Scharlemann~\cite{rub:compar} have shown that the points
$(t, s_t) \in [-1,1] \times [-1,1]$ form a graph $G \subset \mathbf{R}^2$ with 
vertices
\begin{figure}[ht!]\anchor{graphic}\footnotesize
\psfraga <0pt,-2pt> {1}{1}  
\psfraga <-2pt,-4pt> {-1}{${-}1$}  
\psfrag {0}{0}  
\psfrag {f}{$f$}  
\psfraga <-3pt,0pt> {g}{$g$}  
\begin{center}
  \includegraphics[width=1.5in]{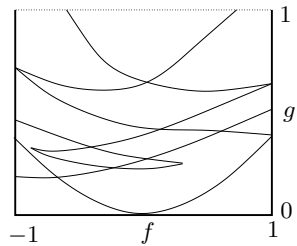}
  \caption{A possible Rubinstein-Scharlemann graphic}
  \label{graphic}
  \end{center}
\end{figure}
of valence $2$ and $4$ in the interior and valence $1$ and $2$ on the edges of
$[-1,1] \times [-1,1]$. We will call this graph the
\textit{Rubinstein-Scharlemann graphic}.

The valence-2 vertices in the interior arise when $g_t$ is a near-Morse function with a
degenerate critical point.
As $t$ passes the vertex, a central singularity and a saddle singularity in
$f_t$ come together and cancel.

At a valence-4 vertex, $(t, s_t)$,
$g_t$ has two critical points at the same level.
As $t$ passes the vertex, one of the critical points rises past the other.
The two critical points may be on the same component of $g_t^{-1}(s_t)$
or different components.

If points $(t,s)$ and $(t', s')$ are in the same component of
$([-1,1] \times [-1,1]) \setminus G$ then there is an isotopy of the surface $f_t$
onto the surface $f_{t'}$ taking the collection of loops $g_t^{-1}(s)$
to the loops $g_{t'}^{-1}(s')$.  In other words, the topological properties
of the loops in the level set depend only on the component $c$, not the
choice of a point within $c$.

In the following section, we will replace the sweep-out $g$ with
a sweep-out of a ball $B \subset M$, i.e.\ a map $g : B \rightarrow [0,1]$
such that $g^{-1}(0)$ is a single point in $B$ and $g^{-1}(s)$ is a sphere
parallel to $\partial B$ for every $s \in (0,1]$.

In this case, the definition of general position is essentially the same.
The function $g_t$ is defined on a closed submanifold of $f_t$. The resulting
graphic, as in \figref{graphic}, looks like the bottom half of the
Rubinstein-Scharlemann graphic for two sweep-outs.

\section{How a sweep-out intersects a closed ball}
\label{lemsect}

Let $(\Sigma, H_1, H_2)$ be a Heegaard splitting of $M$ with spine $K$.
Let $B \subset M$ be a closed ball such that $K$ is transverse to
$\partial B$.  Let $g$ be the genus of $\Sigma$.

\begin{Prop}
\label{mainlem}
Suppose for each surface $f_t$ there is an essential curve of $f_t$
lying completely in $B$. If $g = 1$ then $\Sigma$ is stabilized and $M = S^3$.  If $g \geq 2$ then $\Sigma$ is weakly reducible.
\end{Prop}

\begin{proof}
Choose a point $x \in B$ away from $K$.
Let $N$ be a regular neighborhood of $B$ and
let $g : N \rightarrow [0,1]$ be a sweep out of $N$,  such
that $g^{-1}(0) = \{x\}$.
By work of Cerf~\cite{cerf:diffeo}, we can isotope $g$ by an arbitrarily small
amount so that $g$ and $f$ are in general position.

In particular,
isotope and rescale $g$ so that $B \subset g^{-1}([0,1]) \subset N$.
Let $\{(t_i, s_i)\}$ be the vertices of the Rubinstein-Scharlemann graphic,
ordered so that $t_i < t_{i+1}$.  (i.e.\ $f_{t_i}$ has either a degenerate
critical point at $s_i$ or two non-degenerate critical points at level $s_i$.)

The assumption that there is an essential curve of $f_t$ lying inside $B$
is equivalent to assuming that some level curve of $g_t$ is essential in $f_t$.
If $g_t$ is Morse and each level curve is trivial then $f_t \cap B$ is either a 
collection of disjoint disks or a sphere.  In either case, there are no
essential loops of $f_t$ lying inside $B$.  Thus for the remainder of the 
proof, we will assume that for each $t$, some level curve of $g_t$ is
essential in $f_t$.

Let $G$ be the Rubinstein-Scharlemann graphic in $[-1,1] \times [0,1]$.
For each component $c$ of $([-1,1] \times [0,1]) \setminus G$, let $(t_c, s_c)$ be
a point in the interior of $c$.  Let $L_c$ be the collection of loops
$g_{t_c}^{-1}(s_c)$. If one of the loops bounds a properly embedded, essential
disk in the handlebody $f^{-1}[-1,t_c]$,
label the region $c$ with a $1$ and if a loop bounds a disk in $f^{-1}[t_c,1]$, label
$c$ with a $2$.  This label is independent of our choice of $(t_c, s_c)$
within $c$. Before continuing, we need the following lemma:

\begin{Lem}
Let $(\Sigma, H_1, H_2)$ be a Heegaard splitting for $M$ and let $S$ be a
2-sphere embedded in $M$ such that $S \cap \Sigma$ is a non-empty collection
of loops.  If one or more of the loops is essential in $\Sigma$ then
there is a loop of $S \cap \Sigma$ which bounds a properly embedded, essential
disk in either $H_1$ or $H_2$.
\end{Lem}

\begin{proof}
Induct on the number of components of $S \cap \Sigma$.  If there is a single
component $l$ of $S \cap \Sigma$ then $l$ must be essential.  This $l$
splits $S$ into two disks, each of which is a properly embedded, essential
disk in one of the handlebodies.

If there are more than one components of $S \cap \Sigma$, let $n$ be the
number of loops and assume the result is true whenever there are $n - 1$
loops.  Let $l$ be an innermost loop of $S \cap \Sigma$ (one that bounds
a disk $D$ in $S$ disjoint from $\Sigma$.)  If $l$ is essential, then $D$
is is properly embedded and essential and we are done.

Otherwise, $D$ is boundary parallel (since handlebodies are irreducible) so
we can isotope $S$ to a surface $S'$ such that $S' \cap \Sigma$ is a proper
subset of $S \cap \Sigma$.  Since $S' \cap \Sigma$ has $n - 1$ components,
one of them bounds a properly embedded, essential disk, so the corresponding
component of $S \cap \Sigma$ does as well.
\end{proof}

If follows from the lemma that if a component $c$ of $([-1,1] \times [0,1]) \setminus G$
is not labeled then all the loops of $L_c$ are trivial in $f_t$.  If one of 
the components has both labels $1$ and $2$, then the Heegaard splitting is weakly reducible, so we can assume each component has at most one label.

If a vertical line in $I \times I$ intersects both a component labeled with a 
$1$ and a component labeled with a $2$ then the Heegaard splitting is weakly
reducible:  Choose $t$ so that the vertical line $\{t\} \times I$ intersects 
both types of components.  Then one level set of $g_{t}$ will contain a loop 
bounding a disk in $H_1$ and another level set will contain a loop bounding a disk
in $H_2$.  The level curves are disjoint so the disks will have disjoint
boundaries.

Assume $\Sigma$ is not weakly irreducible.  We will show that $\Sigma$ must
be genus-one and stabilized.
For each $t \in I$, if $g_t$ is a Morse function then there is an essential
level curve of $g_t$. So, each vertical line $\{t\} \times I$ must pass through
either a vertex of $G$ or one or more labeled components, all with the same 
label.  In the latter case, let $l_t$ be this label.

For $t$ near $-1$, $f_t$ approaches the spine $K_1$, so
$l_t$ must be $1$ for small enough $t$. (Recall, $K$ is transverse to
$\partial B$.)  Similarly, $l_t$ must be $2$
for $t$ near $1$.  Thus the label
must switch from $1$ to $2$ somewhere.  It can only change at
a vertex of $G$, so let $v = (t_i, s_i)$ be a vertex where
$l_{t_i - \epsilon}$ is $1$ and $l_{t_i + \epsilon}$ is $2$.

The vertical line at $t_i$ passes through only one vertex.
Because each vertical line to the left or
\begin{figure}[ht!]\anchor{crossing}\small
\psfrag{1}{1}  
\psfrag{2}{2}  
\begin{center}
  \includegraphics[width=1in]{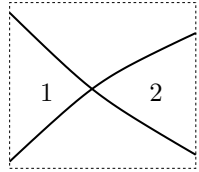}
  \caption{A crossing at which $l_c$ switches, without implying that
    $\Sigma$ is weakly reducible}
  \label{crossing}
  \end{center}
\end{figure}
right of $t_i$ can only pass through one type of component,
the vertex $v$ must be a crossing with the regions above and below it 
unlabeled.
This type of crossing is shown in \figref{crossing}.

\begin{Lem}
\label{toruslem}
Let $g$ be a Morse function on a compact, connected, closed, orientable surface $F$ with critical values $c_i$.
For some $i$, let $h_1 \in (c_{i},c_{i+1})$, $h_2 \in (c_{i+1}, c_{i+2})$
and $h_3 \in (c_{i+2}, c_{i+3})$.  If $g^{-1}(h_1) \cup g^{-1}(h_3)$
consists of trivial loops and $g^{-1}(h_2)$ contains non-trivial loops then
$F$ is a torus.
\end{Lem}

\begin{proof}
The level set $g^{-1}(h_2)$ is the result of adding a band to some loop
$a_1 \subset g^{-1}(h_1)$.  By assumption, $a_1$ must be trivial so let
$D_1 \subset F$ be the disk whose boundary is $a_1$.  Let $D_2 \subset F$
be the 1-handle which is attached to $D_1$.

The level set $g^{-1}(h_3)$ is the result of adding a second 1-handle.
Because all the loops at this level are trivial, the band must be attached
to $\partial (D_1 \cup D_2)$.  Let $D_3$ be this 1-handle and let
$a_3 = \partial (D_1 \cup D_2 \cup D_3)$.  Then $a_3$ must bound a disk $D_4$
because $a_3$ is isotopic to a loop in $g^{-1}(h_3)$.

This $D_4$ must be disjoint from the other handles so the four disks
suggest a handle decomposition of $F$.  Therefore the Euler characteristic
of $F$ is $0$ and $F$ is a torus.
\end{proof}

We can now complete the proof of Proposition~\ref{mainlem}.
If $l_c$ changes at a crossing as in \figref{crossing} then the Morse 
function $g_{t_i - \epsilon}$ satisfies the conditions of Lemma~\ref{toruslem}, implying
$\Sigma$ is a torus. This proves Proposition~\ref{mainlem} when the genus of 
$\Sigma$ is greater than one.
The only case that remains is when the Heegaard splitting has genus one.

Let $g_1 = g_{t_i - \epsilon}$ and let $g_2 = g_{t_i + \epsilon}$.  The only
difference between the two functions is the order of the 2-handles that
create and destroy the essential loops.  One can check that the loops created
by adding the 1-handles in different orders will intersect in a single point
and bound disks in $H_1$, $H_2$, respectively.

Thus if $\Sigma$ is not weakly reducible then $\Sigma$ is a stabilized, genus-one Heegaard splitting (of $S^3$).
\end{proof}

\begin{Coro}[Waldhausen's Theorem]
\label{waldcoro}
Every positive-genus Heegaard splitting of $S^3$ is stabilized.
\end{Coro}

\begin{proof}
Let $(\Sigma, H_1, H_2)$ be a Heegaard splitting of $S^3$ with genus $g$.
Let $f$ be a sweep-out for $\Sigma$.  Let $N$ be a small ball in $S^3$
such that $f^{-1}(t) \cap N$ is either a point, a disk or empty for all $t$.
Then $B = S^3 \setminus N$ is a ball and every level set of $f$ intersects
$B$ non-trivially.  

If $g = 1$ then $\Sigma$ is stabilized and we are done.  If $g \geq 2$ then $\Sigma$ is weakly reducible.
The main theorem in~\cite{cass:red} states that if a Heegaard splitting is weakly reducible then it is reducible or the ambient manifold is Haken.  
Because $S^3$ contains no incompressible surfaces, this implies that $\Sigma$ is reducible. 

Thus $\Sigma$ is a connect sum of Heegaard splittings of $S^3$ with strictly lower genus.  By induction, $\Sigma$ is a connect sum of $g$ splittings, each with genus 1.  Because each of these splittings is stabilized, their connect sum $\Sigma$ is stabilized.
\end{proof}

Note that the proof of Corollary~\ref{waldcoro} is a modification of 
Theorem 5.11 in~\cite{rub:compar}.

\section{Proof of the Main Theorem}
\label{lastsect}

We will begin by stating three lemmas which will be necessary in the final
proof.  A proof of the first lemma will be left to the reader.

\begin{Lem}
\label{induclem}
Let $T$ be a graph tangle in a ball, $B$, and let $D \subset B$ be a properly
embedded disk, disjoint from $T$, cutting $B$ into two balls $B_1$ and $B_2$.
Assume that $T \cap B_1$ is unknotted in $B_1$.  Then $T$ is unknotted in $B$
if and only if $T \cap B_2$ is unknotted in $B_2$.
\end{Lem}

Let $\Sigma$ be a strongly irreducible Heegaard splitting and let
$K$ and $B$ be as in Section~\ref{lemsect}. Assume that $K$ has no trivial 
cones with respect to $B$.

\begin{Lem}
\label{planlem}
If $\Sigma$ is strongly irreducible and $B \subset M$ is a closed ball
then there is an immersion of $K \times I$ such that $K \times [0,1)$ is
embedded, $K \times \{0\} = K$ and $K \times \{1\}$ is disjoint from $B$.
\end{Lem}

\begin{proof}
Let $f$ be a sweep-out for $\Sigma$.  By Proposition~\ref{mainlem}, there must be a
value $t \in (-1,1)$ such that $B \cap f_t$ does not contain any essential
loops of $f_t$.  Let $L$ be the collection
of loops $f_t \cap \partial B$.  These loops are pairwise disjoint and each
loop $\alpha \in L$ bounds a disk $D_\alpha \subset f_t$.

There is an embedding of $K_1 \times [0,1]$ in $f^{-1}[-1,t]$ such that
$K_1 \times \{0\} = K_1$ and
$K_1 \times \{1\} \subset f_t \setminus \bigcup D_\alpha$.
Then $K_1 \times \{1\}$ is disjoint from $L$ so $K_1 \times \{1\}$
is either in $B$ or disjoint from $B$.  If it were inside $B$, then $K_1$
would be isotopic into $B$ (a contradiction) so it must be disjoint from
$B$.  The same construction for $K_2$ finishes the proof.

Note that in general, $K_1 \times \{1\} \cap K_2 \times \{1\}$ will be
non-empty.  The statement of the lemma requires only that the interior of
$K \times I$ be embedded.
\end{proof}
 
Let $K \times [0,1]$ be immersed in $M$ with embedded interior such that 
$K \times \{0\} = K$ and $K \times \{1\}$ is disjoint from $B$, as in
Lemma~\ref{planlem}.  We will write $K \times I$ for the image of this immersion.  Let $L = (K \times I) \cap \partial B$. 

Let $e$ be an edge of $K$.
Then $L \cap (e \times I)$ consists of a number of arcs 
and loops. Because 
$K \times \{1\}$ is disjoint from $B$, each arc in $e \times I$ must have both its 
endpoints in $K \times \{0\}$ or in $v \times (0,1)$ where $v$ is an endpoint
of $e$. Let $\alpha$ be such an arc.

If $\alpha$ lies in $e \times [0,1]$ and has both its endpoints on the same
$v \times (0,1)$ for some $v$, then the endpoints of $\alpha$ bound a sub-arc $\beta$ of $v \times (0,1)$.  The loop $\alpha \cup \beta$ bounds a disk 
$D \subset e \times [0,1]$.  
Sliding $\beta$ across $D$ gives an isotopy of $v \times (0,1)$ which extends to $K \times (0,1)$.   This isotopy reduces the number
of intersections of $v \times [0,1]$ with $\partial B$.  Thus by a finite
number of such isotopies, we can eliminate all arcs with both endpoints on the same vertical edge $v \times (0,1)$.

\begin{Lem}
\label{ltimeslem}
Let $l$ be a component of $L$ which intersects $K \times \{0\}$.
Then $l$ separates $K \times I$ into a component homeomorphic to
$K \times I$ and a component homeomorphic to $l \times I$.
\end{Lem}

\begin{proof}
Consider a projection $\rho$ of $L$ into $K \times \{0\}$.
After eliminating arcs in $L$ with both endpoints on the same edge
$v \times I$, we can choose a projection $\rho$ so that $\rho$ is 
one-to-one on each arc of $L$.
We will prove that $\rho$ is then one-to-one on the whole component $l$.
Assume for contradiction this map sends points $x,y \in l$ to the
same point in $K \times \{0\}$ and assume $y$ is below $x$ (i.e.\ farther from 
$K \times \{0\}$.)

Since $l$ is connected, there is a path $\alpha: [0,1] \rightarrow l$ 
with $\alpha(0) = x$ and $\alpha(1) = y$.  We would like to construct
a path $\alpha' : [0,1] \rightarrow l$ such that 
$\rho \circ \alpha (t) = \rho \circ \alpha' (t)$ for every $t$
and $\alpha'(0) = y$.  This is, in some sense, a lift of the projection
of $\alpha$, starting at a different base point.  We will show that
for each $t$, $\alpha'(t)$ is below $\alpha(t)$.

The edges of $v \times I$ of $K \times I$ cut $\alpha$ into a number of arcs coinciding with the arcs of $L$. The first arc 
can't end on $K \times \{0\}$ because then $\alpha$ would end. It can't
end at $y$ because $\rho$ is one-to-one on each arc of $l$.  Thus 
the first arc must pass from $x$ to some $v \times I$.  

The point $y$ also sits on an arc of $l$. This arc extends in the same 
direction as the first arc of $\alpha$.  It can't end on $K \times \{0\}$
because it is below another arc of $L$ which extends to $v \times I$. It
cannot end on $K \times \{1\}$ because $L$ is disjoint from $K \times \{1\}$.
Thus the arc ends on $v \times I$, below the first arc of $\alpha$.
We can define $\alpha'$ as desired, strictly below the first arc of $\alpha$.

Continue in this fashion for each arc in the image of $\alpha$, until
$\alpha$ reaches $y$.  Now, $\alpha'(1)$ is a point $z$ below $y$ such that 
$\rho(z) = \rho(y) = \rho(x)$.  Continuing the construction, we can find
another point below $z$, and so on, implying that $\rho^{-1} \rho (x)$ 
contains an infinite number of points. But this is impossible because $L$
is a finite graph.  Thus $\rho$ must be one-to-one.

For each point $x \in l$, consider the vertical arc $a_x$ from $x$
to $K \times \{0\}$.  Because $\rho$ is one-to-one, the union of
these arcs is homeomorphic to $l \times I$.  Its complement is
homeomorphic to $K \times I$.
\end{proof}

We are now ready to prove the main theorem.  Isotope $K \times (0,1)$
so as to minimize the number of loops in $L$ and the number of
intersections with vertical edges $v \times I$.  Each component of $K
\cap B$ contains one or more components of $L$.  If there are no
loops, then each $l \times I$ component lies entirely in $B$ so by
Lemma~\ref{ltimeslem}, each $l$ cuts off a piece of $K \times I$
homeomorphic $l \times I$.  This defines an isotopy of a component of
$K \cap B$ into $\partial B$ so $K \cap B$ is unknotted.  Otherwise,
if there are loops in $L$, we will show that there is always a disk
properly embedded in $B$ which cuts off an unknotted tangle, allowing
us to use Lemma~\ref{induclem}.

Let $\beta$ be a loop of $L$.  Then $\beta$ bounds a disk, $D_1 \subset (K \times [0,1])$.  Assume $\beta$ is an inner-most loop, i.e.\ assume $D_1
\cap \partial B = \partial D_1$.  There are two cases to consider:
either $D_1$ lies entirely in $B$ or the interior of $D_1$ lies
outside of $B$.

Assume the interior of $D_1$ is disjoint from $B$.  The loop $\beta$ separates
$\partial B$ into two disks.  Because $M$ allows a strongly irreducible
Heegaard splitting, $M$ is irreducible, so $D_1 \cup \partial B$ separates
$M$ into two balls and a punctured $M$.  Let $D_2$ be the disk of
$\partial B \setminus \beta$ such that $D_1 \cup D_2$ bounds a ball, $B'$
outside $B$.

The graph $K \times \{1\}$ is disjoint from $\partial B'$ so each
component of $K \times \{1\}$ must lie outside of $B'$. (If it
were inside $B'$ then $K \times I$ would suggest an isotopy
bringing one of the spines into $B'$.) 

If $K \times I$ intersects $D_2$, then let $l$ be a component of $D_2 \cap K \times I \subset L$.  By Lemma~\ref{ltimeslem}, $l$ cuts off a piece of $K \times I$ homeomorphic to $l \times I$.  This piece is outside of $B$ so $l \times I$ defines
a trivial cone. Thus there are no arc-intersections and $K \times
\{0\}$ lies outside of $B'$.  Any components of $K \times I$ in $B'$ can be
isotoped out of $B'$ and then $D_1$ can be isotoped into $B$,
reducing the number of components of $L$.

So, if $L$ is minimized, we can assume no loops bound disks whose
interiors are disjoint from $B$. Let $\mathbf{D}$ be the
collection of innermost disks in $L$. These separate $B$ into a
number of balls. Let $D \in \mathbf{D}$ be an outermost disk in
$B$, cutting off a ball $B'$ which does not contain any disks of
$\mathbf{D}$. Let $B''$ be the complement of $B'$.  Note that a trivial cone of 
$K$ with respect to $B'$ or $B''$ is also a trivial cone with respect to $B$.

If there are no arcs of $K$ in $B'$, then we can isotope all the
pieces of $K \times I$ out of $B'$ and then isotope off $D$.  Otherwise,
let $L' = (K \times I) \cap \partial B'$ and let $\beta$ be an
innermost loop of $L'$ cutting off a disk $D'$.

The disk $D'$ cannot lie in $B'$ because $D$ was assumed to be an outermost
disk.  Thus $D'$ must lie outside $B'$.  By applying the argument above, to
$B'$ instead of $B$, it is clear that $D'$ can be isotoped into $B'$
without affecting $K$.  So, if the number of loops of intersection is
minimized, there can be no loops in $L'$.

Thus $K \cap B'$ is an unknotted tangle in $B'$.  By Lemma~\ref{induclem}, $K \cap B$
is unknotted if and only if $K \cap B''$ is unknotted.
Again minimize the number of loops of intersection for  this new tangle.
The number of components of $(K \times I) \cap (B \setminus B')$ is strictly less than
the number of components of $(K \times I) \cap B$ so the process will eventually
terminate, proving the theorem.
\endproof

\Addresses\recd

\end{document}